\documentclass[12pt,thmsa]{article}
\usepackage{sw20lart}
\usepackage{latexsym}

\input tcilatex

\begin{document}

\title{The largest linear space of operators satisfying the Daugavet Equation in $%
L_1$}
\author{R. V. Shvidkoy \\
Department of Mathematics\\
University of Missouri - Columbia\\
Columbia, MO 65211\\
USA\\
\textit{E-mail:} mathgr31@showme.missouri.edu}
\date{March 14, 1999}
\maketitle

\begin{abstract}
We find the largest linear space of bounded linear operators on $L_1(\Omega
) $ that being restricted to any $L_1(A)$, $A\subset \Omega $, satisfy the
Daugavet equation.
\end{abstract}

\section{Introduction.}

Let $(\Omega ,\Sigma ,\mu )$ be an arbitrary measure space without atoms of
infinite measure. Let also $\Sigma ^{+}=\{A\in \Sigma :\mu (A)>0\}$. If $%
A\in \Sigma ^{+}$, $L_1(A)$ stands for the space of (classes of) $\mu $%
-integrable functions supported on $A$. If $T$ is a bounded linear operator
on $L_1(\Omega )$ and $A\in \Sigma ^{+}$, we denote by $T_A$ the restriction
of $T$ onto $L_1(A)$. Finally, $\mathcal{L}(L_1(\Omega ))$ denotes the space
of all bounded linear operators on $L_1(\Omega ).$

The purpose of this note is to give an explicit description of the largest
linear space $\mathcal{M}$ of operators $T\in \mathcal{L}(L_1(\Omega ))$
satisfying the following identity: 
\begin{equation}
\Vert Id_A+T_A\Vert =1+\Vert T_A\Vert ,  \label{DE}
\end{equation}
for any set $A\in $ $\Sigma ^{+}$.

Identity (\ref{DE}) is known as the Daugavet equation and is investigated in
a series of works (see \cite{KSSW} and \cite{Shv} for recent results and
further references). It was first discovered by Babenko and Pichugov (\cite
{Bab-Pich}) that all the compact operators on $L_1[0,1]$ satisfy (\ref{DE}),
if $A=[0,1]$. Later, Holub proved the same result for the weakly compact
operators on an arbitrary atomless $L_1(\Omega )$ (see \cite{Hol}). Plichko
and Popov in their work \cite{Plich-Pop} found much broader (in case of
atomless $\mu $) linear class of so-called narrow operators satisfying the
Daugavet equation, and in fact their proof works for operators from $L_1(A)$
to $L_1(\Omega )$, whenever $A\in \Sigma ^{+}$.

So, finding the largest class of such operators naturally completes this
line of results.

\section{Main result.}

In the sequel it is convenient to denote $\Sigma _A^{+}=\{B:B\subset A,$ $%
B\in \Sigma ^{+}\}$, whenever $A\in \Sigma ^{+}$.

We define $\mathcal{M}$ as the set of all operators $T\in \mathcal{L}%
(L_1(\Omega ))$ that meet the following condition: 
\begin{eqnarray}
&&\text{for every }\varepsilon >0\text{ and }A\in \Sigma ^{+}\text{ there is
a }B\in \Sigma _A^{+}\text{ with }\mu (B)<\infty \text{ such that}
\label{def} \\
&&\left\| \chi _B\cdot T\left( \frac{\chi _B}{\mu (B)}\right) \right\|
<\varepsilon \text{.}  \nonumber
\end{eqnarray}
This condition simply means that the operator $T$ can shift sufficiently
many functions from their supports.

Let us state our main result.

\begin{theorem}
\label{main}Every linear set of operators satisfying (\ref{DE}) for any $%
A\in \Sigma ^{+}$ is contained in $\mathcal{M}$, and $\mathcal{M}$ is itself
a closed linear space consisting of such operators.
\end{theorem}

The main ingredient in the proof of this theorem is the following
proposition.

\begin{proposition}
\label{aux}For an operator $T\in \mathcal{L}(L_1(\Omega ))$ the following
conditions are equivalent:

\begin{enumerate}
\item[(i)]  $T$ and $-T$ satisfy (\ref{DE}) for all $A\in \Sigma ^{+}$;

\item[(ii)]  For every $\varepsilon >0$ and $A\in \Sigma ^{+}$ there is an $%
A^{\prime }\in \Sigma _A^{+}$ such that if $B\in \Sigma _{A^{\prime }}^{+}$
then we can find a $B^{\prime }\in \Sigma _B^{+}$ with the following
properties:

a) $\left\| \frac{\chi _{B^{\prime }}}{\mu (B^{\prime })}-\frac{\chi _B}{\mu
(B)}\right\| <\varepsilon $,

b) $\left\| \chi _{B^{\prime }}\cdot T\left( \frac{\chi _{B^{\prime }}}{\mu
(B^{\prime })}\right) \right\| <\varepsilon $;

\item[(iii)]  $T\in \mathcal{M}$.
\end{enumerate}
\end{proposition}

\emph{Proof.} (i) implies (ii). We begin with the following observation.

Suppose $S:L_1(A)\mapsto L_1(\Omega )$ is a bounded linear operator, then
any given $\varepsilon >0$ there is a set $A_1\in \Sigma _A^{+}$ with $\mu
(A_1)<\infty $ such that for every non-negative function $f\in S(L_1(A_1))$
we have $\Vert Sf\Vert >\Vert S\Vert -\varepsilon $.

Indeed, we can assume that $\mu (A)<\infty $ and choose $g^{*}\in
S(L_1^{*}(\Omega ))$ so that $\Vert S^{*}g^{*}\Vert >\Vert S\Vert
-\varepsilon $. Then, regarding $S^{*}g^{*}$ as an element of $L_\infty (A)$
we find a set $A_1\in \Sigma _A^{+}$ with $\theta S^{*}g^{*}(A_1)\subset
(\Vert S\Vert -\varepsilon ,\Vert S\Vert ]$, where $\theta $ is a sign. Now,
if $f\in S(L_1(A))$, $f\geq 0$ and $\limfunc{supp}(f)\subset A_1$, then $%
\Vert Sf\Vert >\theta g^{*}(Sf)=\theta S^{*}g^{*}(f)>\Vert S\Vert
-\varepsilon $, from where the observation follows.

We know that $\Vert Id_A+T_A\Vert =1+\Vert T_A\Vert $. By scaling, without
loss of generality we can and do assume that $\Vert T_A\Vert =1$. So there
is an $A_1\in \Sigma _A^{+}$ with $\mu (A_1)<\infty $ such that 
\begin{equation}
\left\| \frac{\chi _B}{\mu (B)}+T\left( \frac{\chi _B}{\mu (B)}\right)
\right\| >2-\varepsilon \text{,}  \label{daug1}
\end{equation}
whenever $B\in \Sigma _{A_1}^{+}$. We also know that $\Vert
Id_{A_1}-T_{A_1}\Vert =1+\Vert T_{A_1}\Vert >2-\varepsilon $. Thus there
exists an $A^{\prime }\in \Sigma _{A_1}^{+}$ such that 
\begin{equation}
\left\| \frac{\chi _B}{\mu (B)}-T\left( \frac{\chi _B}{\mu (B)}\right)
\right\| >2-\varepsilon \text{,}  \label{daug2}
\end{equation}
whenever $B\in \Sigma _{A^{\prime }}^{+}$.

We prove that $A^{\prime }$ is the desired set.

To this end, let us fix $B\in \Sigma _{A^{\prime }}^{+}$. It follows from (%
\ref{daug1}), (\ref{daug2}) and a theorem of Dor \cite{Dor} that there are
two disjoint measurable sets $\Omega _1$ and $\Omega _2$ in $\Omega $ such
that 
\begin{equation}
\int_{\Omega _1}\left| T\left( \frac{\chi _B}{\mu (B)}\right) \right|
(t)dt>(1-\varepsilon )^2\text{,}  \label{eq3}
\end{equation}
and 
\[
\int_{\Omega _2}\frac{\chi _B}{\mu (B)}(t)dt>(1-\varepsilon )^2\text{.} 
\label{eq4} 
\]
The last inequality implies 
\begin{eqnarray}
\mu (B\cap \Omega _1) &=&\mu (B)\int_{B\cap \Omega _1}\frac{\chi _B}{\mu (B)}%
(t)dt<\mu (B)\int_{\Omega \backslash \Omega _2}\frac{\chi _B}{\mu (B)}(t)dt
\label{eq5} \\
\ &<&(1-(1-\varepsilon )^2)\mu (B)=(2\varepsilon -\varepsilon ^2)\mu (B)%
\text{.}  \nonumber
\end{eqnarray}
Let us put $B^{\prime }=B\backslash \Omega _1$ and show that $B^{\prime }$
meets conditions a) and b).

First, 
\begin{eqnarray*}
\left\| \frac{\chi _{B^{\prime }}}{\mu (B^{\prime })}-\frac{\chi _B}{\mu (B)}%
\right\| =\int_\Omega \left| \frac{\chi _{B^{\prime }}}{\mu (B^{\prime })}-%
\frac{\chi _{B^{\prime }}}{\mu (B)}+\frac{\chi _{B^{\prime }}}{\mu (B)}-%
\frac{\chi _B}{\mu (B)}\right| (t)dt && \\
\leq 1-\frac{\mu (B^{\prime })}{\mu (B)}+\frac{\mu (B\cap \Omega _1)}{\mu (B)%
}=2\frac{\mu (B\cap \Omega _1)}{\mu (B)} &&\text{,}
\end{eqnarray*}
and taking into account (\ref{eq5}), we obtain 
\begin{equation}
\left\| \frac{\chi _{B^{\prime }}}{\mu (B^{\prime })}-\frac{\chi _B}{\mu (B)}%
\right\| <2(2\varepsilon -\varepsilon ^2)\text{.}  \label{eq20}
\end{equation}
Second, from (\ref{eq3}), (\ref{eq20}) and $\Vert T_A\Vert =1$ it follows
that 
\begin{eqnarray*}
\left\| \chi _{B^{\prime }}\cdot T\left( \frac{\chi _{B^{\prime }}}{\mu
(B^{\prime })}\right) \right\| &=&\int_{B^{\prime }}\left| T\left( \frac{%
\chi _{B^{\prime }}}{\mu (B^{\prime })}\right) \right| (t)dt \\
\ &<&\int_{B^{\prime }}\left| T\left( \frac{\chi _B}{\mu (B)}\right) \right|
(t)dt+2(2\varepsilon -\varepsilon ^2) \\
\ &\leq &\int_{\Omega \backslash \Omega _1}\left| T\left( \frac{\chi _B}{\mu
(B)}\right) \right| (t)dt+2(2\varepsilon -\varepsilon ^2) \\
\ &\leq &3(2\varepsilon -\varepsilon ^2)\text{.}
\end{eqnarray*}
In view of arbitrariness of $\varepsilon $, this gives the desired result.

It is obvious that (iii) follows from (ii).

Let us finally prove that (iii) implies (i). Since $\mathcal{M}$ is stable
under scalar multiplication, it is sufficient to prove (\ref{DE}) only for $%
T $.

To this end, we fix an arbitrary $A\in \Sigma ^{+}$ and as above for any
given $\varepsilon >0$ we find an $A^{\prime }\in \Sigma _A^{+}$ with $\mu
(A^{\prime })<\infty $ such that for every $B\in \Sigma _{A^{\prime }}^{+}$, 
$\left\| T\left( \frac{\chi _B}{\mu (B)}\right) \right\| >\Vert T_A\Vert
-\varepsilon $. By condition (\ref{def}), there is a $B_0\in \Sigma
_{A^{\prime }}^{+}$ such that $\left\| \chi _{B_0}\cdot T\left( \frac{\chi
_{B_0}}{\mu (B_0)}\right) \right\| <\varepsilon $. This means that $\frac{%
\chi _{B_0}}{\mu (B_0)}$ and $T\left( \frac{\chi _{B_0}}{\mu (B_0)}\right) $
are almost disjoint functions, and as a consequence we have the following
estimate:

\begin{eqnarray*}
\Vert Id_A+T_A\Vert &\geq &\left\| \frac{\chi _{B_0}}{\mu (B_0)}+T\left( 
\frac{\chi _{B_0}}{\mu (B_0)}\right) \right\| \\
&=&\int_{B_0}\left| \frac{\chi _{B_0}}{\mu (B_0)}+T\left( \frac{\chi _{B_0}}{%
\mu (B_0)}\right) \right| (t)dt+\int_\Omega \left| T\left( \frac{\chi _{B_0}%
}{\mu (B_0)}\right) \right| (t)dt \\
&&-\int_{B_0}\left| T\left( \frac{\chi _{B_0}}{\mu (B_0)}\right) \right|
(t)dt \\
&>&1-\varepsilon +\Vert T_A\Vert -\varepsilon -\varepsilon =1+\Vert T_A\Vert
-3\varepsilon \text{.}
\end{eqnarray*}
This finishes the proof. $\Box $

\bigskip Now we are in a position to prove our main result.

\smallskip\emph{Proof of Theorem \ref{main}.}

Proposition \ref{aux} easily implies that $\mathcal{M}$ is largest and
consists of operators satisfying (\ref{DE}) for all $A\in \Sigma ^{+}$. $%
\mathcal{M}$ is obviously closed and stable under scaling. So, the only
thing we have to prove is that if operators $U$ and $V$ belong to $\mathcal{M%
}$, then their sum belong to $\mathcal{M}$ too. To show this, we check
condition (ii) of Proposition \ref{aux} for $U+V$. Further on, we assume
that $\Vert V\Vert \leq 1$.

Indeed, let $A\in \Sigma ^{+}$ and $\varepsilon >0$ be arbitrary. Applying
Proposition \ref{aux} to the operator $U$ we find a set $A^{\prime }\in
\Sigma _A^{+}$ as in condition (ii). Then, by the same proposition applied
to $V$ we find a set $A^{\prime \prime }\in \Sigma _{A^{\prime }}^{+}$ with
the correspondent properties. To show that $A^{\prime \prime }$ is the
required set, suppose $B\in \Sigma _{A^{\prime \prime }}^{+}$. By the choice
of $A^{\prime \prime }$ there is a $B^{\prime }\in \Sigma _B^{+}$ such that 
\begin{equation}
\left\| \frac{\chi _{B^{\prime }}}{\mu (B^{\prime })}-\frac{\chi _B}{\mu (B)}%
\right\| <\frac \varepsilon 4,  \label{eq10}
\end{equation}
and 
\begin{equation}
\left\| \chi _{B^{\prime }}\cdot V\left( \frac{\chi _{B^{\prime }}}{\mu
(B^{\prime })}\right) \right\| <\frac \varepsilon 4.  \label{eq12}
\end{equation}
Since $B^{\prime }\subset A^{\prime }$, by the analogous property of $%
A^{\prime }$, there is a $B^{\prime \prime }\in \Sigma _{B^{\prime }}^{+}$
with 
\begin{equation}
\left\| \frac{\chi _{B^{\prime \prime }}}{\mu (B^{\prime \prime })}-\frac{%
\chi _{B^{\prime }}}{\mu (B^{\prime })}\right\| <\frac \varepsilon 4,
\label{eq11}
\end{equation}
and 
\[
\left\| \chi _{B^{\prime \prime }}\cdot U\left( \frac{\chi _{B^{\prime
\prime }}}{\mu (B^{\prime \prime })}\right) \right\| <\frac \varepsilon 2. 
\]
From (\ref{eq10}) and (\ref{eq11}) we get $\left\| \frac{\chi _{B^{\prime
\prime }}}{\mu (B^{\prime \prime })}-\frac{\chi _B}{\mu (B)}\right\|
<\varepsilon $. So, if we prove that 
\[
\left\| \chi _{B^{\prime \prime }}\cdot V\left( \frac{\chi _{B^{\prime
\prime }}}{\mu (B^{\prime \prime })}\right) \right\| <\frac \varepsilon 2%
\text{,} 
\]
then 
\[
\left\| \chi _{B^{\prime \prime }}\cdot (V+U)\left( \frac{\chi _{B^{\prime
\prime }}}{\mu (B^{\prime \prime })}\right) \right\| <\varepsilon \text{,} 
\]
and we are done. But this easily follows from (\ref{eq12}), (\ref{eq11}) and
the facts that $\Vert V\Vert \leq 1$ and $B^{\prime \prime }\subset
B^{\prime }$.

The proof is completed. $\Box $

\end{document}